\documentclass{amsart}
\usepackage{amsfonts}
\usepackage{amssymb}
\usepackage{xypic}

\newcommand{\Z}{{\mathbb Z}}
\newcommand{\C}{{\mathbb C}}
\renewcommand{\P}{{\mathbb P}}
\newcommand{\QH}{\text{\sl QH}}

\newcommand{\Gr}{\operatorname{Gr}}

\renewcommand{\mod}{\operatorname{mod}}
\newcommand{\bull}{{\sssize \bullet}}
\renewcommand{\O}{{\mathcal O}}

\newtheorem{lemma}{Lemma} 
 
\newtheorem{thm}{Theorem} 
\newtheorem{cor}{Corollary}
\newcommand{\reflemma}[1]{Lemma~\ref{#1}}
\newcommand{\refthm}[1]{Theorem~\ref{#1}}
\newcommand{\refcor}[1]{Corollary~\ref{#1}}

\newcommand{\rmc}{\hat}
\newcommand{\rmr}{\bar}
\newcommand{\rmrc}[1]{\rmc{\rmr #1}}

\begin{document}

\title{Quantum cohomology of Grassmannians}
\author{Anders Skovsted Buch}
\address{Massachusetts Institute of Technology \\
  Building 2, Room 248 \\
  77 Massachusetts Avenue \\
  Cambridge, MA 02139
}
\date{\today}
\email{abuch@math.mit.edu}
\thanks{The author was partially supported by NSF Grant DMS-0070479}
\maketitle

\section{Introduction}

The purpose of this paper is to give simple proofs of the main
theorems about the (small) quantum cohomology ring of a Grassmann
variety.  This first of all includes Bertram's quantum versions of the
Pieri and Giambelli formulas \cite{bertram:quantum}.  Bertram's proofs
of these theorems required the use of quot schemes.  Our proof of the
quantum Pieri formula uses no moduli spaces and only the definition of
Gromov-Witten invariants.  In fact we show that this formula is a
consequence of the classical Pieri formula.  We then show that the
quantum Giambelli formula follows immediately from the quantum Pieri
formula together with the classical Giambelli formula and
associativity of quantum cohomology \cite{ruan.tian:mathematical,
kontsevich.manin:gromov-witten}.

We also give a short proof of the Grassmannian case of a formula of
Fulton and Woodward for the minimal $q$-power which appears in a
quantum product of two Schubert classes \cite{fulton.woodward:fixme}.
In addition we supply a proof of a simple version of the rim-hook
algorithm \cite{bertram.ciocan-fontanine.ea:quantum} based on ``mod
$n$'' arithmetic which is due to F.~Sottile \cite{sottile:rational}.
Finally we recover Siebert and Tian's presentation of the quantum
cohomology of Grassmannians \cite{siebert.tian:on*1}.

In this paper we only assume associativity of quantum cohomology and
standard facts about the usual cohomology.  The basic idea is that if
a rational curve of degree $d$ passes through a Schubert variety in
the Grassmannian $\Gr(l,\C^n)$ of $l$-dimensional subspaces of $\C^n$,
then the linear span of the points of this curve gives rise to a point
in $\Gr(l+d,\C^n)$ which must lie in a related Schubert variety.
Remarkably, this simple idea can in many cases be used to conclude
that no curves pass through three Schubert varieties in general
position because the intersection of the related Schubert varieties in
$\Gr(l+d,\C^n)$ is empty.  In particular, the quantum Giambelli
formula can be deduced by knowing that certain Gromov-Witten
invariants are zero, and in each case this follows because the
codimensions of the related Schubert varieties add up to more than the
dimension of $\Gr(l+d,\C^n)$.  In a paper
\cite{buch.kresch.ea:quantum} with A.~Kresch and H.~Tamvakis the same
idea will be applied to obtain a similar treatment of the quantum
cohomology of Lagrangian and orthogonal Grassmannians
\cite{kresch.tamvakis:qclg, kresch.tamvakis:qort}.

We thank Sottile for triggering our search for simplifications to the
theory of quantum cohomology for Grassmannians, when he explained his
simplified rim-hook algorithm to us.  We also thank W.~Fulton for
numerous helpful comments to our paper, some of which greatly enhanced
the clarity of our proof of the quantum Giambelli formula.  Finally we
thank Kresch for simplifying our proof of the key
\reflemma{lemma:span}.

\section{Preliminaries}


Set $E = \C^n$, $X = \Gr(l,E)$, and $k = n-l$.  Given a flag of
subspaces $F_1 \subset F_2 \subset \dots \subset F_n = E$ and a
partition $\lambda = (\lambda_1 \geq \lambda_2 \geq \dots \geq
\lambda_l \geq 0)$ with $\lambda_1 \leq k$, we define the Schubert
variety
\begin{equation} \label{eqn:schubert}
  \Omega_\lambda(F_\bull) = \{ V \in X \mid 
   \dim(V \cap F_{k+i-\lambda_i}) \geq i ~\forall 1 \leq i \leq l \} \,.
\end{equation}
The codimension of this variety is equal to the weight $|\lambda| =
\sum \lambda_i$ of $\lambda$.  We let $\Omega_\lambda$ denote the
class of $\Omega_\lambda(F_\bull)$ in the cohomology ring $H^*X =
H^*(X; \Z)$.  The Schubert classes $\Omega_\lambda$ form a basis for
this ring.  The partitions indexing this basis are exactly those where
the Young diagram fits in an $l \times k$ rectangle.

The cohomology ring has the presentation $H^*X = \Z[\Omega_1, \dots,
\Omega_k]/(Y_{l+1}, \dots, Y_n)$ where $Y_p = \det(\Omega_{1+j-i})_{1
  \leq i,j \leq p}$.  Here we set $\Omega_i = 0$ for $i < 0$ or $i >
k$ for convenience.  In this presentation the class $\Omega_\lambda$
is given by the Giambelli formula
\begin{equation} \label{eqn:giambelli}
  \Omega_\lambda = \det(\Omega_{\lambda_i+j-i})_{1\leq i,j \leq l} \,.
\end{equation}
This is usually deduced ({cf.\ }\cite{fulton:young}) from the Pieri
formula, which states that
\begin{equation} \label{eqn:pieri}
  \Omega_i \cdot \Omega_\lambda = \sum \Omega_\nu
\end{equation}
where the sum is over all partitions $\nu$ which can be obtained by
adding $i$ boxes to the Young diagram of $\lambda$ with no two in the
same column.


Recall that the degree of a rational curve $f : \P^1 \to \P^N$ is
equal to the number of points in the inverse image by $f$ of a general
hyperplane in $\P^N$.  The degree of a map $f : \P^1 \to X$ is the
degree of the composition of $f$ with the Pl{\"u}cker embedding $X
\subset \P(\bigwedge^l E)$ which maps a point $V \in X$ to $v_1 \wedge
\dots \wedge v_l$ for any basis $\{v_1, \dots, v_l\}$ of $V$.  Now let
$\lambda$, $\mu$, and $\nu$ be three partitions contained in an $l
\times k$ rectangle, and let $d \geq 0$ be an integer such that
$|\lambda| + |\mu| + |\nu| = l k + d n$.  The Gromov-Witten invariant
$\langle \Omega_\lambda, \Omega_\mu, \Omega_\nu \rangle_d$ is defined
as the number of rational curves of degree $d$ on $X$, which meet all
of the Schubert varieties $\Omega_\lambda(F_\bull)$,
$\Omega_\mu(G_\bull)$, and $\Omega_\nu(H_\bull)$ for general flags
$F_\bull$, $G_\bull$, $H_\bull$, up to automorphisms of $\P^1$.  A
simple proof that this number is well defined is given in
\cite{bertram:quantum}.  If $|\lambda| + |\mu| + |\nu| \neq l k + d n$
then $\langle \Omega_\lambda, \Omega_\mu, \Omega_\nu \rangle_d = 0$.

The (small) quantum cohomology ring $\QH^*X = \QH^*(X;\Z)$ of $X$ is a
$\Z[q]$-algebra which is isomorphic to $H^*X \otimes_\Z \Z[q]$ as a
module over $\Z[q]$.  In this ring we have Schubert classes
$\sigma_\lambda = \Omega_\lambda \otimes 1$.  The ring structure on
$\QH^*X$ is defined by
\begin{equation} \label{eqn:qstructure}
  \sigma_\lambda \cdot \sigma_\mu = \sum_{\nu,\, d \geq 0}
  \langle \Omega_\lambda, \Omega_\mu, \Omega_{\nu^\vee} \rangle_d \, q^d \,
  \sigma_\nu
\end{equation}
where $\nu^\vee = (k-\nu_l, k-\nu_{l-1}, \dots, k-\nu_1)$ is the
partition for the dual Schubert class of $\Omega_\nu$.  It is a
non-trivial fact that this defines an associative ring structure
\cite{ruan.tian:mathematical, kontsevich.manin:gromov-witten} (see
also \cite{fulton.pandharipande:notes}).  Notice that the definition
implies that $\QH^*X$ is a graded ring where $\sigma_\lambda$ has
degree $|\lambda|$ and $q$ has degree $n$.  Furthermore, the map
$\QH^*X/(q) \to H^*X$ which sends $\sigma_\lambda$ to $\Omega_\lambda$
is an isomorphism of rings, so the quantum ring is a deformation of
the usual cohomology ring.

Siebert and Tian have given a presentation of this ring which is
similar to the above presentation of the cohomology ring
\cite{siebert.tian:on*1}.  We will comment on this presentation
towards the end of this paper.  For now we will start from the basic
assumption that this quantum ring is a well defined associative ring,
and aim towards proving generalizations of the Pieri and Giambelli
formulas stated above to reveal its structure.

\section{The span and kernel of a curve}

Our main new tool is the following definition.  If $Y$ is any
subvariety of $X = \Gr(l,E)$ we define the {\em span\/} of $Y$ to be
the smallest subspace of $E$ containing all the $l$-dimensional spaces
given by points of $Y$.  Similarly we define the {\em kernel\/} of $Y$
to be the largest subspace of $E$ contained in all the spaces given by
points of $Y$.

\begin{lemma} \label{lemma:span}
  Let $C$ be a rational curve of degree $d$ in $X$.  Then the span of
  $C$ has dimension at most $l+d$ and the kernel of $C$ has dimension
  at least $l-d$.
\end{lemma}

\begin{proof}
  Let $C$ be the image of a regular function $f : \P^1 \to X$ of
  degree $d$, and let $S \subset E \otimes \O_X$ be the tautological
  subbundle on $X$.  Then $f^* S = \bigoplus_{i=1}^l \O_{\P^1}(-a_i)$
  for integers $a_i \geq 0$ with sum $d$, and $f$ is given by an
  inclusion $\bigoplus_{i=1}^l \O_{\P^1}(-a_i) \subset E \otimes
  \O_{\P^1}$, i.e.\ a point $p \in \P^1$ is mapped to the fiber over
  $p$ of the image of this bundle map.  If $(s:t)$ are homogeneous
  coordinates on $\P^1$ then $\Gamma(\O_{\P^1}(a_i))$ has the basis
  $\{ s^j t^{a_i-j} \}_{0\leq j\leq a_i}$, so each map
  $\O_{\P^1}(-a_i) \to E \otimes \O_{\P^1}$ has the form
\[
  \sum_{j=0}^{a_i} \alpha_j \, s^{-j} t^{j-a_i} ~\mapsto~ 
  \sum_{j=0}^{a_i} v^{(i)}_j \otimes \alpha_j
\]
for vectors $v^{(i)}_j \in E$ (which will depend on the chosen
identification of $f^*S$ with $\bigoplus_{i=1}^l \O_{\P^1}(-a_i)\,$).
The span of $C$ must therefore be contained in the span of the set $\{
v^{(i)}_j \}$ which has cardinality $\sum 1+a_i = l+d$.  On the other
hand, at least $l-d$ of the integers $a_i$ must be zero, and the
kernel of $C$ contains the span of the corresponding vectors
$v^{(i)}_0$.
\end{proof}

The above lemma can also be obtained by proving that any regular
function $f : \P^1 \to X$ can be written as $f(t) = f_1(t) \wedge
\cdots \wedge f_l(t)$ for maps $f_i : \P^1 \to \P(E)$.  This is what
we did until Kresch showed us the simpler argument given here.


If $\lambda$ is a partition and $d$ a non-negative integer, we let
$\rmc \lambda$ denote the partition obtained by removing the leftmost
$d$ columns from the Young diagram of $\lambda$, {i.e.\ }$\rmc
\lambda_i = \max(\lambda_i-d,0)$.

\begin{lemma}
Let $C \subset X$ be a rational curve of degree $d \leq k$ and let $W
\subset E$ be an $l+d$ dimensional subspace containing the span of
$C$.  If $\lambda$ is a partition such that $C \cap
\Omega_\lambda(F_\bull) \neq \emptyset$ then $W$ belongs to the
Schubert variety $\Omega_{\rmc
\lambda}(F_\bull)$ in $\Gr(l+d, E)$.
\end{lemma}
\begin{proof}
Let $V \in C \cap \Omega_\lambda(F_\bull)$.  Since $V \subset W$, the
Schubert conditions on $V$ imply that $\dim(W \cap F_{k+i-\lambda_i})
\geq i$ for all $i$, which says exactly that $W$ belongs to
$\Omega_{\rmc \lambda}(F_\bull)$.
\end{proof}

Similarly, if $S \subset E$ is an $l-d$ dimensional subspace contained
in the kernel of $C$, then $S \in \Omega_{\rmr \lambda}(F_\bull)
\subset \Gr(l-d, E)$, where $\rmr \lambda = (\lambda_{d+1}, \dots,
\lambda_l)$ is the result of removing the top $d$ rows of $\lambda$.
Therefore the condition that a curve meets a given set of Schubert
varieties in $X$ implies that intersections of related Schubert
varieties in $\Gr(l+d,E)$ and in $\Gr(l-d,E)$ are not empty, which is
a statement about the usual cohomology of these spaces.  As we shall
see, this simple idea is sufficient to compute Gromov-Witten
invariants in many important cases.


\section{The quantum Pieri formula}

We start with the following quantum version of the Pieri formula
\cite{bertram:quantum}.

\begin{thm}[Bertram] \label{thm:qpieri}
If $\lambda$ is contained in an $l \times k$ rectangle and $p \leq k$ then
\[ \sigma_p \cdot \sigma_\lambda = \sum \sigma_\mu + q \sum \sigma_\nu \]
where the first sum is over all partitions $\mu$ such that $|\mu| =
|\lambda|+p$ and $k \geq \mu_1 \geq \lambda_1 \geq \mu_2 \geq
\lambda_2 \geq \dots \geq \mu_l \geq \lambda_l$, and the second sum is
over all partitions $\nu$ such that $|\nu| = |\lambda| + p - n$ and
$\lambda_1 - 1 \geq \nu_1 \geq \lambda_2 - 1 \geq \nu_2 \geq \dots
\geq \lambda_l - 1 \geq \nu_l \geq 0$.
\end{thm}

Recall that the length $\ell(\lambda)$ of a partition $\lambda$ is the
number of non-zero parts of $\lambda$.  Notice that the second sum in
the theorem is non-zero only if $\ell(\lambda) = l$.

\begin{proof}
  The first sum is dictated by the classical Pieri formula.  Notice
  that this classical case is equivalent to the following statement.
  If $\alpha$ and $\beta$ are partitions such that $|\alpha| + |\beta|
  + p = l k$ then
\[ \langle \Omega_\alpha, \Omega_\beta, \Omega_p \rangle_0 =
   \begin{cases} 
     1 & \text{if $\alpha_i + \beta_j \geq k$ for $i+j =
       l$ and $\alpha_i + \beta_j \leq k$ for $i+j = l+1$;} \\
     0 & \text{otherwise.}
   \end{cases}
\]

Now suppose $|\alpha| + |\beta| + p = l k + d n$ for some $d \geq 1$
and let $C$ be a rational curve of degree $d$ in $X$ which meets each
of the varieties $\Omega_\alpha(F_\bull)$, $\Omega_\beta(G_\bull)$,
and $\Omega_p(H_\bull)$ for general flags $F_\bull$, $G_\bull$,
$H_\bull$.  Let $W \subset E$ be a subspace of dimension $l+d$ which
contains the span of $C$.  Then $W \in \Gr(l+d, E)$ lies in the
intersection $\Omega_{\rmc \alpha}(F_\bull) \cap \Omega_{\rmc
\beta}(G_\bull) \cap \Omega_{\rmc p}(H_\bull)$ where $\rmc \alpha$ and
$\rmc \beta$ are the results of removing the leftmost $d$ columns from
$\alpha$ and $\beta$, and $\rmc p = \max(p-d, 0)$.  Since the flags
$F_\bull$, $G_\bull$, $H_\bull$ are general, this implies that $|\rmc
\alpha| + |\rmc \beta| + \rmc p \leq (l+d)(k-d)$.  Since we also have
\[ |\rmc \alpha| + |\rmc \beta| + \rmc p
   \geq |\alpha| + |\beta| - 2ld + p - d = (l+d)(k-d) + d^2 - d
\]
we deduce that $d = 1$ and $\ell(\alpha) = \ell(\beta) = l$.  Using
this, the quantum Pieri formula becomes equivalent to the statement
that if $|\alpha| + |\beta| + p = l k + n$ then $\langle \Omega_\alpha,
\Omega_\beta, \Omega_p \rangle_1 = 1$ if $\alpha_i + \beta_j \geq k+1$
for $i+j = l+1$ and $\alpha_i + \beta_j \leq k+1$ for $i+j = l+2$;
otherwise $\langle \Omega_\alpha, \Omega_\beta, \Omega_p \rangle_1 =
0$.  In other words, the quantum Pieri formula states that $\langle
\Omega_\alpha, \Omega_\beta, \Omega_p \rangle_1 = \langle \Omega_{\rmc
  \alpha}, \Omega_{\rmc \beta}, \Omega_{\rmc p} \rangle_0$ where the
right hand side is a coefficient of the classical Pieri formula for
$\Gr(l+1, E)$.

If $\langle \Omega_{\rmc \alpha}, \Omega_{\rmc \beta}, \Omega_{\rmc p}
\rangle_0$ is zero then the space $W$ can't exist, so neither can $C$.
On the other hand, if $\langle \Omega_{\rmc \alpha}, \Omega_{\rmc
    \beta}, \Omega_{\rmc p} \rangle_0 = 1$ then there exists a unique
subspace $W \subset E$ of dimension $l+1$ which is contained in the
intersection $\Omega_{\rmc \alpha}(F_\bull) \cap \Omega_{\rmc
  \beta}(G_\bull) \cap \Omega_{\rmc p}(H_\bull)$.  Furthermore, since
the flags are general, $W$ must lie in the interior of each of these
Schubert varieties.  In particular, each of the spaces $V_1 = W \cap
F_{n - \alpha_l}$ and $V_2 = W \cap G_{n - \beta_l}$ have dimension
$l$.  Notice also that $V_1 \in \Omega_\alpha(F_\bull)$ and $V_2 \in
\Omega_\beta(G_\bull)$.  Since $\Omega_\alpha(F_\bull) \cap
\Omega_\beta(G_\bull) = \emptyset$ we deduce that $V_1 \neq V_2$, so
$S = V_1 \cap V_2$ has dimension $l-1$.  We conclude that the only
rational curve of degree $1$ in $X$ which meets the Schubert varieties
for $\alpha$, $\beta$, and $p$ is the line $\P(W/S)$ of
$l$-dimensional subspaces between $S$ and $W$.
\end{proof}

\section{The quantum Giambelli formula}

Now set $\sigma_i = 0$ for $i < 0$ and for $k < i < n$.  The quantum
Giambelli formula states that the classical formula
(\ref{eqn:giambelli}) continues to be valid if all Schubert classes
are replaced with the corresponding quantum Schubert classes
\cite{bertram:quantum}.

\begin{thm}[Bertram] \label{thm:qgiambelli}
  If $\lambda$ is a partition contained in an $l \times k$ rectangle,
  then the Schubert class $\sigma_\lambda$ in $\QH^* \Gr(l,E)$ is
  given by $\sigma_\lambda = \det(\sigma_{\lambda_i+j-i})_{1 \leq i,j
    \leq l}$.
\end{thm}
\begin{proof}
  We claim that if $0 \leq i_j \leq k$ for $1 \leq j \leq l$ then
  $\sigma_{i_1} \cdot \sigma_{i_2} \cdots \sigma_{i_l} = (\Omega_{i_1}
  \cdot \Omega_{i_2} \cdots \Omega_{i_l}) \otimes 1$, i.e.\ no
  $q$-terms show up when the first product is expanded in the quantum
  ring.  Using induction, this can be established by proving that if
  $\ell(\mu) < l$ then the expansion of $\sigma_i \cdot \sigma_\mu$
  involves no $q$-terms and no partitions of lengths greater than
  $\ell(\mu)+1$.  The claim therefore follows from
  \refthm{thm:qpieri}.  Since the determinant of the quantum Giambelli
  formula is a signed sum of products of the form $\sigma_{i_1} \cdot
  \sigma_{i_2} \cdots \sigma_{i_l}$, we conclude from the classical
  Giambelli formula (\ref{eqn:giambelli}) that
  $\det(\sigma_{\lambda_i+j-i}) = \det(\Omega_{\lambda_i+j-i}) \otimes
  1 = \Omega_\lambda \otimes 1 = \sigma_\lambda$ as required.
\end{proof}

Notice that this proof uses only that no $q$-terms are contained in a
product $\sigma_i \cdot \sigma_\mu$ when $\ell(\mu) < l$ (in addition
to the classical Giambelli and Pieri formulas).  In fact,
\refthm{thm:qpieri} could be replaced with the following lemma.

\begin{lemma} \label{lemma:noq}
  Let $\lambda$ and $\mu$ be partitions contained in an $l \times k$
  rectangle such that $\ell(\lambda) + \ell(\mu) \leq l$.  Then
  $\sigma_\lambda \cdot \sigma_\mu = (\Omega_\lambda \cdot \Omega_\mu)
  \otimes 1$.
\end{lemma}
\begin{proof}
  If $d \geq 1$ and $\nu$ is a partition such that $|\lambda| + |\mu|
  + |\nu| = l k + nd$, then any intersection $\Omega_{\rmc
    \lambda}(F_\bull) \cap \Omega_{\rmc \mu}(G_\bull) \cap
  \Omega_{\rmc \nu}(H_\bull)$ of general Schubert varieties in
  $\Gr(l+d,E)$ must be empty since $|\rmc \lambda| + |\rmc \mu| +
  |\rmc \nu| \geq |\lambda| + |\mu| + |\nu| - 2 d l = l k + d k - d l >
  (l+d)(k-d)$.  This shows that $\langle \Omega_\lambda, \Omega_\mu,
    \Omega_\nu \rangle_d = 0$.
\end{proof}

A dual version of this lemma states that if $\lambda_1 + \mu_1 \leq k$
then $\sigma_\lambda \cdot \sigma_\mu = (\Omega_\lambda \cdot
\Omega_\mu) \otimes 1$.  This follows by replacing the span of a curve
with its kernel in the above proof, or by using the duality
isomorphism $\QH^* \Gr(l,E) \cong \QH^* \Gr(k,E^*)$.

\section{The minimal power of $q$ in a quantum product}

As a further demonstration of the use of the span of a rational curve
in $X$, we will give a short proof of a recent theorem of Fulton and
Woodward \cite{fulton.woodward:fixme} in the case of Grassmannians.
It generalizes the fact that any product $\sigma_\lambda \cdot
\sigma_\mu$ in $\QH^*X$ is non-zero.

Given a partition $\lambda$ and an integer $d \geq 0$, we let
$\rmrc{\lambda}$ denote the partition obtained by removing the top $d$
rows and the leftmost $d$ columns from $\lambda$.  Recall that a
product $\Omega_\lambda \cdot \Omega_\mu$ is non-zero in $H^*
\Gr(l,E)$ if and only if $\lambda_i + \mu_j \leq k$ for $i+j =
l+1$.

\begin{thm}[Fulton and Woodward]
  The smallest power of $q$ which appears in a product $\sigma_\lambda
  \cdot \sigma_\mu$ in $\QH^*X$ is equal to the smallest $d$ for
  which $\Omega_{\rmrc{\lambda}} \cdot \Omega_{\mu} \neq 0$ in $H^*X$.
\end{thm}
\begin{proof}
  If $\sigma_\lambda \cdot \sigma_\mu$ contains $q^d$ times a Schubert
  class then some curve of degree $d$ meets each of the Schubert
  varieties $\Omega_\lambda(F_\bull)$ and $\Omega_\mu(G_\bull)$ where
  $F_\bull$ and $G_\bull$ are general flags.  If $W \subset E$ has
  dimension $l+d$ and contains the span of this curve then $W$ lies in
  the intersection $\Omega_{\rmc \lambda}(F_\bull) \cap \Omega_{\rmc
    \mu}(G_\bull)$ in $\Gr(l+d,E)$.  In particular this intersection
  is not empty which implies that $\Omega_{\rmc \lambda} \cdot
  \Omega_{\rmc \mu} \neq 0$ in $H^* \Gr(l+d,E)$.  Since this is
  equivalent to $\Omega_{\rmrc{\lambda}} \cdot \Omega_\mu \neq 0$ in
  $H^*X$, this proves the inequality ``$\geq$'' of the theorem.
  
  Now let $d$ be the smallest number for which
  $\Omega_{\rmrc{\lambda}} \cdot \Omega_\mu \neq 0$.  Notice that this
  implies that $\lambda$ contains a $d \times d$ rectangle, {i.e.\ 
    }$\lambda_d \geq d$.  Set $\alpha = (k+d-\lambda_d, \dots,
  k+d-\lambda_1)$ and let $\beta$ be the partition given by $\beta_i =
  \max(d-\lambda_{l+1-i}, 0)$ for $1 \leq i \leq l$.  If the Young
  diagram for $\lambda$ is put in the upper-left corner of an $l$ by
  $k+d$ rectangle, then $\alpha$ is the complement of $\lambda$ in the
  top $d$ rows of this rectangle, turned $180$ degrees, and $\beta$ is
  the complement of $\lambda$ in the leftmost $d$ columns, also
  turned.
  
  It follows from the Littlewood-Richardson rule that the product
  $\sigma_\lambda \cdot \sigma_\beta$ contains the class
  $\sigma_{(d^l)+\rmc \lambda} = \sigma_{(d^l)} \cdot \sigma_{\rmc
  \lambda}$ and that $\sigma_{\rmc \lambda} \cdot \sigma_\alpha$
  contains $\sigma_{(k^d),\rmrc{\lambda}} = \sigma_{(k^d)} \cdot
  \sigma_{\rmrc{\lambda}}$.  Since the structure constants of $\QH^*
  X$ are all non-negative, this implies that $\sigma_\lambda \cdot
  \sigma_\beta \cdot \sigma_\alpha$ contains the product
  $\sigma_{(d^l)} \cdot \sigma_{(k^d)} \cdot \sigma_{\rmrc{\lambda}} =
  q^d \sigma_{\rmrc{\lambda}}$.  Since $\Omega_{\rmrc{\lambda}} \cdot
  \Omega_\mu \neq 0$ by assumption, we conclude that $\sigma_\lambda
  \cdot \sigma_\mu \cdot \sigma_\alpha \cdot \sigma_\beta$ contains
  $q^d$ times some Schubert class.  In particular, at least one term
  of the product $\sigma_\lambda \cdot \sigma_\mu$ must involve a
  power of $q$ which is less than or equal to $d$.  This proves the
  other inequality ``$\leq$'' in the theorem.  Notice that the
  identities we have used in this argument follow easily from
  \refthm{thm:qpieri} and the dual version of \reflemma{lemma:noq},
  combined with the classical Pieri rule.
\end{proof}

\section{The rim-hook algorithm}

We will next recall how to carry out computations in the quantum
cohomology ring of $X$.  Let $c_i \in \QH^*X$ denote the $i$th Chern
class of the dual of the tautological subbundle $S \subset E \otimes
\O_X$.  This means that $c_i = \sigma_{(1^i)}$ for $0 \leq i \leq l$
and $c_i = 0$ for $i < 0$ and $i > l$.  For $p \geq 1$ we then define
$\sigma_p = \det(c_{1+j-i})_{1 \leq i,j \leq p}$ in $\QH^*X$.  Notice
that for $p < n$ we have $\sigma_p = \Omega_p \otimes 1$, so this
definition is compatible with our previous definition of $\sigma_p$.
The definition implies that for every $m \geq 1$ we have
\begin{equation} \label{eqn:sigmap}
\sum_{i=0}^l (-1)^i\, \sigma_{m-i} \, \sigma_{(1^i)} = 0 \,. 
\end{equation}

\begin{lemma} \label{lemma:shift}
  For any $p \geq k+1$ we have $\sigma_p = (-1)^{l-1}\, q\,
  \sigma_{p-n}$.
\end{lemma}
\begin{proof}
  Both sides of the identity are zero for $k+1 \leq p < n$.  It
  follows from (\ref{eqn:sigmap}) with $m = n$ that $\sigma_n + (-1)^l
  \sigma_k\, \sigma_{(1^l)} = 0$, so $\sigma_n = (-1)^{l-1} \sigma_k\,
  \sigma_{(1^l)} = (-1)^{l-1} q$ by \refthm{thm:qpieri}.  For $p > n$
  we finally obtain
\[
\sigma_p 
= \sum_{i=1}^l (-1)^{i-1} \sigma_{p-i}\, \sigma_{(1^i)} 
= (-1)^{l-1}\, q\, \sum_{i=1}^l (-1)^{i-1} \sigma_{p-n-i}\, \sigma_{(1^i)}
= (-1)^{l-1}\, q\, \sigma_{p-n}
\]
by induction, which proves the lemma.
\end{proof}

If $I = (I_1, I_2, \dots, I_p)$ is any sequence of integers we set
$\sigma_I = \det(\sigma_{I_i+j-i})_{1 \leq i,j \leq p}$ in $\QH^*X$.
This is compatible with our notation $\sigma_\lambda$ for partitions
$\lambda$ by \refthm{thm:qgiambelli}.  Any element $\sigma_I$ can be
rewritten as zero or plus or minus $\sigma_\mu$ for some partition
$\mu$ by moves of the form $\sigma_{I,a,a+1,J} = 0$ and
$\sigma_{I,a,b,J} = - \sigma_{I,b-1,a+1,J}$, which correspond to
interchanging rows in the determinant defining $\sigma_I$.  Notice
that if $\lambda$ is a partition then we have the formal identity
$\sigma_\lambda = \det(c_{\lambda'_i+j-i})$ where $\lambda'$ is the
conjugate partition of $\lambda$, obtained by mirroring the Young
diagram for $\lambda$ in its diagonal.  In particular, if
$\ell(\lambda) > l$ then the top row of the latter determinant is
zero, and $\sigma_\lambda = 0$.

\begin{cor} \label{cor:rimhook}
  Let $\lambda$ be a partition with at most $l$ parts.  For each $1
  \leq j \leq l$, choose $I_j \in \Z$ such that $I_j \equiv \lambda_j
  ~(\mod n)$ and $j-l \leq I_j < j-l+n$.  Then we have
\[ \sigma_\lambda = (-1)^{d(l-1)}\, q^d \, \sigma_I  \in \QH^*X \]
where $I = (I_1, \dots, I_l)$ and $n d = |\lambda| - \sum I_j$.
\end{cor}


Notice that the moves described above will rewrite the chosen element
$\sigma_I$ to zero or plus or minus a basis element $\sigma_\mu$ of
$\QH^*X$ with $\mu$ contained in an $l \times k$ rectangle.  The
corollary is equivalent to the dual rim-hook algorithm of
\cite{bertram.ciocan-fontanine.ea:quantum}.  The formulation in terms
of reduction modulo $n$ is due to F.~Sottile \cite{sottile:rational},
who in turn was inspired by the work of Ravi, Rosenthal, and Wang
\cite{ravi.rosenthal.ea:degree}.  This formulation gives an efficient
algorithm to determine if two elements $\sigma_\lambda$ and
$\sigma_\mu$ are identical.  In particular it shows that
$\sigma_\lambda = 0$ if and only if $\lambda_i - i \equiv \lambda_j -
j ~(\mod n)$ for some $i \neq j$.

\refcor{cor:rimhook} makes it easy to do computations in $\QH^*X$
(cf.\ \cite{goodman.wenzl:littlewood-richardson, walton:fusion,
  bertram.ciocan-fontanine.ea:quantum} and
\cite[Ex.~13.35]{kac:infinite-dimensional*6}).  In order to expand a
product $\sigma_\lambda \cdot \sigma_\mu$, one uses the classical
Littlewood-Richardson rule to write $\sigma_\lambda \cdot \sigma_\mu =
\sum c^\nu_{\lambda \mu} \sigma_\nu$ where the sum is over partitions
$\nu$ of length at most $l$ and $c^\nu_{\lambda \mu}$ is the
Littlewood-Richardson coefficient.  Each of the elements $\sigma_\nu$
can then be reduced to a basis element of $\QH^*X$ by the corollary.

\section{Generators and relations}

Finally, we show how to deduce Siebert and Tian's presentation of the
quantum cohomology ring of $X$ \cite{siebert.tian:on*1}.  We claim
that
\begin{equation} \label{eqn:qpresent}
  \QH^*X = \Z[c_1, \dots, c_l, q]/
  (\sigma_{k+1}, \dots, \sigma_{n-1}, \sigma_n+(-1)^l q) \,.
\end{equation}
In fact, all of the given relations vanish in $\QH^*X$ by
\reflemma{lemma:shift}.  On the other hand, the proof of
\reflemma{lemma:shift} together with \refcor{cor:rimhook} shows that
the linear map $H^*X \otimes \Z[q] \to
\Z[c_1,\dots,c_l,q]/(\sigma_{k+1}, \dots, \sigma_{n-1},
\sigma_n+(-1)^l q)$ which sends $\Omega_\lambda \otimes q^d$ to $q^d
\sigma_\lambda$ is surjective, which proves the isomorphism.  Siebert
and Tian gave this presentation in its dual form $\QH^*X =
\Z[\sigma_1, \dots, \sigma_k, q]/(\tilde Y_{l+1}, \dots, \tilde
Y_{n-1}, \tilde Y_n + (-1)^k q)$ where the $\tilde Y_i$ are defined
inductively by $\tilde Y_0 = 1$, $\tilde Y_i = 0$ for $i < 0$, and
$\tilde Y_i = \sum_{j=1}^k (-1)^{j-1} \sigma_j \tilde Y_{i-j}$ for $i
> 0$.

The presentation (\ref{eqn:qpresent}) can also be stated as $\QH^*X =
\Z[c_1,\dots,c_l]/I^{(l,n)}$ where $I^{(l,n)} = (\sigma_{k+1}, \dots,
\sigma_{n-1})$.  In this form the variable $q$ is represented by
$(-1)^{l-1} \sigma_n = \sigma_{(k+1,1^{l-1})}$.  In
\cite{goodman.wenzl:littlewood-richardson} other generators for the
ideal $I^{(l,n)}$ are used, namely all elements $\sigma_\lambda$ for
which $\lambda_1 - \lambda_l = k+1$.  It is an easy exercise to show
that these elements in fact generate the same ideal, and to show that
the basis $\{ q^d \sigma_\lambda \mid \lambda_1 \leq k \}$ of $\QH^*X$
is equal to the set $\{ \sigma_\lambda \mid \lambda_1 - \lambda_l \leq
k \}$ which is the basis used in
\cite{goodman.wenzl:littlewood-richardson}.  We refer to
\cite{sottile:rational} for a more detailed discussion of the
relations of quantum cohomology of Grassmannians with other fields.


\providecommand{\bysame}{\leavevmode\hbox to3em{\hrulefill}\thinspace}
\providecommand{\href}[2]{#2}

\end{document}